\documentclass[10pt]{article}
\usepackage{amssymb}
\usepackage{amsmath,amssymb}\usepackage{cite}\usepackage{mathrsfs}\usepackage[amsmath,thmmarks]{ntheorem}
\usepackage{listings}
\usepackage[titletoc]{appendix}\allowdisplaybreaks

\makeatletter
\renewcommand\theequation{\thesection.\@arabic\c@equation}
\makeatother

\newtheorem{thm}{ Theorem}[section]%
\newtheorem{lem}[thm]{ Lemma}%
\newtheorem{Exam}[thm]{ Example}%
\newtheorem{Def}[thm]{Definition}%
\newtheorem{Que}[thm]{Question}%
\newtheorem{Con}[thm]{Conjecture}%
\input cyracc.def

\def\f{\noindent}

\def\demo{\f{\bf Proof}\hskip10pt}

\def\qed{\hfill $\Box$}

\begin{document}
\title{\textbf{The generating pairs of the 2-transitive groups}}
\author{Junyao Pan\,
 \\\\
School of Mathematical and Statistics, Huaiyin Normal University, Huai'an, Jiangsu,\\ 223300
 People's Republic of China \\}
\date {} \maketitle

\baselineskip=16pt

\vskip0.5cm

{\bf Abstract:} Given a finite group $G$. The generating pair $(H,a)$ of $G$, that is, $H<G$ and $a\in G$ such that $\langle a,H\rangle=G$. In this paper, we introduce the definition of FF-subgroup to characterize the generating pairs of the symmetric groups, alternating groups and projective groups $PSL(2,q)$. This
gives a partial answer to an open problem of J. Andr\'e and J. Ara\'ujo and P. J. Cameron.



\section {Introduction}

It is well-known that the generating sets for groups are far more complicated than generating sets for vector spaces. So it has always been an interesting subject to study the generating sets of groups. Recently, J. Andr\'e and J. Ara\'ujo and P. J. Cameron proposed the following open problem about the generating sets for groups in \cite{AAC}.

\begin{Que}\label{pan1-1}\normalfont(\cite[Problem~1]{AAC}~)
 Let $G\leq S_n$ be a 2-transitive group. Classify the generating pairs $(a,H)$, where $a\in S_n$ and $H\leq S_n$, such that $\langle a,H\rangle=G$.
\end{Que}

In fact, the maximal subgroups are closely related to generating sets, and the research of this direction has attracted the attention of some scholars, such as the references \cite{CLR,LS}. Thereby, we also try to solve the Question\ \ref{pan1-1} by applying the maximal subgroups. Note that all elements of Frattini subgroup are not generating elements, and thus $H$ is far from the Frattini subgroup if $(a,H)$ is a generating pair. On the other hand, it is evident that for each maximal subgroup $M$, the $(a,M)$ is a generating pair where $a\notin M$. Motivated by these two features, we introduce the following notions and notations.

Given a proper subgroup $H$ of a group $G$. We denote by $\Delta_H(G)$ the union of all maximal subgroups of $G$ containing $H$, and further we say $\Delta_H(G)$ is the \emph{maximal cover} of $H$ in $G$.

\begin{Def}\label{pan1-2}\normalfont
 Let $H$ be a proper subgroup of a finite group $G$. If $\Delta_H(G)\neq G$, then we say $H$ is a \emph{FF-subgroup} of $G$.
 \end{Def}

 Obviously, every maximal subgroup of a finite group is a FF-subgroup, however, the Frattini subgroup and its subgroups of a finite non-cyclic group are not FF-subgroups. Most critically, the FF-subgroup has the following property.

\begin{lem}\label{pan1-3}\normalfont
Let $G$ be a finite group with a proper subgroup $H$. Then $\langle H,a\rangle=G$ if and only if $a\in G\setminus\Delta_H(G)$ and $H$ is a FF-subgroup of $G$.
\end{lem}
\demo It is simple to show the sufficiency holds. Assume $\langle H,a\rangle=G$ and $H$ is not a FF-subgroup of $G$. Then $\Delta_H(G)=G$, and thus there exists
a maximal subgroup $M$ in $G$ such that $H\leq M$ and $a\in M$, a contradiction. Similarly, we obtain $a\in G\setminus\Delta_H(G)$. The proof of this lemma is completed.   \qed

According to Definition\ \ref{pan1-2} and Lemma\ \ref{pan1-3}, we see FF-subgroups are not only closely related to maximal subgroups but also to generating pairs, and further the Question\ \ref{pan1-1} is equivalent to classify the FF-subgroups of the 2-transitive groups in some sense.

Let's go back to see if the Definition\ \ref{pan1-2} is trivial. In other words, are there some groups which have some subgroups are neither the subgroups of their Frattini subgroups nor FF-subgroups? Actually, there are many such groups, and then we give an example.

\begin{Exam}\label{pan1-4}\normalfont(\cite{J1}~)
Consider the 3-generator group
$G=\langle x, y, z|y^{-1}xy=y^{b-2}x^{-1}y^{b+2}, z^{-1}yz=z^{c-2}y^{-1}z^{c+2}, x^{-1}zx=x^{a-2}z^{-1}x^{a+2} \rangle$
where $a, b, c$ are even integers distinct from zero. Then the subgroups $\langle x\rangle$, $\langle y\rangle$, $\langle z\rangle$ are neither contained in the Frattini subgroup nor FF-subgroups.
\end{Exam}
\demo Since $G$ is a 3-generator group, it follows that there dot not exists an element $g\in G$ such that $\langle g,\langle x\rangle\rangle=G$ or $\langle g,\langle y\rangle\rangle=G$ or $\langle g,\langle z\rangle\rangle=G$. Then by Lemma\ \ref{pan1-3} we see the subgroups $\langle x\rangle$, $\langle y\rangle$, $\langle z\rangle$ are not FF-subgroups. On the other hand, it is easy to see the maximal subgroup containing $y$ and $z$ does not contain $x$, and thus $\langle x\rangle$ is not the subgroup of Frattini subgroup of $G$. Similarly, $\langle y\rangle$ and $\langle z\rangle$ are also not the subgroups of Frattini subgroup of $G$. The proof of this example is completed.   \qed

 In this paper we will show that all nontrivial subgroups of the finite symmetric and alternating groups and projective groups $PSL(2,q)$ are FF-subgroups.

\section {Symmetric and Alternating Groups }

In this section, a permutation has cycle-type $(d_1,d_2,...,d_t)$ where the $d_i$ are distinct, that is, the permutation is the product of $t$ disjoint cycles whose lengths are $d_i$ for $i=1,2,...,t$. Moreover, $S_n$ is the symmetric group acts on $[n]$, where $[n]=\{1,2,...,n\}$.

In \cite{LPS}, M. W. Liebecka, C. E. Praeger and J. Saxl shaw that if $X$ is $A_n$ or $S_n$, acting on a set $\Omega$ of size $n$, and $G$ is any
maximal subgroup of $X$ with $G\neq A_n$, then $G$ satisfies one of the following:

(a) $G=(S_m\times S_{k})\cap X$, with $n=m+k$ and $m\neq k$ (intransitive case);

(b) $G=(S_k wr S_m)\cap X$, with $n=km$, $m>1$ and $k>1$ (imprimitive case);

(c) $G=AGL_k(p)\cap X$, with $n=p^k$ and $p$ is a prime (affine case);

(d) $G=(T^k\cdot(Out T\times S_k))\cap X$, with $T$ a nonabelian simple group, $k\geq2$ and $n=|T|^{k-1}$(diagonal case);

(e) $G=(S_k wr S_m)\cap X$, with $n=k^m$, $k\geq5$ and $m\geq2$, excluding the case where $X=A_n$ and $G$ is imprimitive on $\Omega$;

(f) $T\lhd G\leq Aut(T)$, with $T$ a nonabelian simple group, $T\neq A_n$ and $G$ acting primitively on $\Omega$ (almost simple case).

Given a subgroup $H$ of $S_n$. Then we use $Max(S_n, H)$ to denote the set of all maximal subgroups of $S_n$ containing $H$, see \cite{MG}. Obviously, $H$ is the subgroup of the intersection of all the maximal subgroups in $Max(S_n, H)$, and thus if the intersection of some maximal subgroups in $Max(S_n, H)$ is the identity then $H$ is the identity. Let's first deal with the symmetric groups along this line.

\begin{lem}\label{pan2-1}\normalfont
Let $n$ be an odd number. Then all nontrivial subgroups of $S_n$ are FF-subgroups.
\end{lem}
\demo It is clear that there is no nontrivial subgroup in $S_1$, and all nontrivial subgroups of $S_3$ are maximal subgroups. Hence, it suffices to consider the case $n>3$. Furthermore, we know that the maximal subgroup $S_m\times S_{k}$ contains the element has cycle-type $(n-2,2)$ if and only if $m=2$ or $k=2$. In other words, the element has cycle-type $(n-2,2)$ is contained in a maximal subgroup $S_{n-2}\times S_2$. We claim that every element has cycle-type $(n-2,2)$ is only contained in one maximal subgroup $S_{n-2}\times S_2$.

Since $n$ is an odd number, we see $n-2$ and $2$ are coprime, therefore, if a maximal subgroup contains an element has cycle-type $(n-2,2)$, then the maximal subgroup contains a $2$-cycle. Applying \cite[Corollary 1.3]{J}, it follows that the primitive maximal subgroup contains an element has cycle-type $(n-2,2)$ is $A_n$, however, the element has cycle-type $(n-2,2)$ is an odd permutation. Hence, all primitive maximal subgroups do not contain the element has cycle-type $(n-2,2)$. Considering the imprimitive maximal subgroup $S_k wr S_m$ with $n=mk$, $m>1$ and $k>1$. Then there exists a partition of $[n]$ into $m$ sets of size $k$ such that $S_k wr S_m$ is the stabiliser of this partition. Noticing that $S_k wr S_m$ contains the element has cycle-type $(n-2,2)$ if and only if $k=2$ or $m=2$, and this is contradict to the $n$ is an odd number. Therefore, our claim holds.

Assume that the nontrivial subgroup $H$ of $S_n$ is not a FF-subgroup. Then by Definition\ \ref{pan1-2} we see $\Delta_H(S_n)=S_n$. On the other hand, we note that for each element has cycle-type $(n-2,2)$, there exists unique maximal subgroup $S_{n-2}\times S_2$ contains it. According to our claim, it follows that all the maximal subgroups $S_{n-2}\times S_2$ are contained in $Max(S_n,H)$. Obviously, the intersection of all the maximal subgroups $S_{n-2}\times S_2$ is the identity, and thus $H$ is the identity, a contradiction. We have thus proved this lemma.   \qed

\begin{lem}\label{pan2-2}\normalfont
Let $n$ be an even number. Then all nontrivial subgroups of $S_n$ are FF-subgroups.
\end{lem}
\demo Obviously, it suffices to prove the case $n\geq4$. Since $n$ is an even number, we see $n-3$ and $2$ are coprime. An argument similar to the one used in proving Lemma\ \ref{pan2-1} shows that the element has cycle-type $(n-3,2,1)$ is only contained in the maximal subgroup $S_{n-1}$ or $S_3\times S_{n-3}$ or $S_2\times S_{n-2}$. Note that if $n=4$, then the element has cycle-type $(n-3,2,1)$ is only contained in the maximal subgroup $S_3$, and the proof of Lemma\ \ref{pan2-1} indicates all nontrivial subgroups of $S_4$ are FF-subgroups.

Consider the case $n>4$. Assume that the nontrivial subgroup $H$ of $S_n$ is not a FF-subgroup. However, we observe that $H$ is the identity, in other words, all $x\in[n]$ are fixed points of $H$. Now we start to confirm our observation.

Note that if $x\in[n]$ is a fixed point of a maximal subgroup in $Max(S_n,H)$, then $x$ is the fixed point of $H$. So we suppose $x$ is not the fixed point of all maximal subgroups in $Max(S_n,H)$. Then we see the element with cycle-type $(n-3,2,1)$ and the fixed point $x$ is contained in the maximal subgroup $S_3\times S_{n-3}$ or $S_2\times S_{n-2}$. If all elements with cycle-type $(n-3,2,1)$ and the fixed point $x$ are contained in the maximal subgroups $S_2\times S_{n-2}$, then all these maximal subgroups are contained in $Max(S_n,H)$ and further the intersection of these maximal subgroups is the identity, and so $x$ is the fixed point of $H$. Similarly, $x$ is the fixed point of $H$ if all elements with cycle-type $(n-3,2,1)$ and the fixed point $x$ are contained in the maximal subgroups $S_3\times S_{n-3}$. Now we may assume $\alpha$ with cycle-type $(n-3,2,1)$ and the fixed point $x$ is contained in the maximal subgroup $S_{\Psi}\times S_{[n]\setminus\Psi}$ where $\Psi=\{x,y,z\}$. Since $\Delta_H(S_n)=S_n$, there exist two maximal subgroups $S_{\Delta}\times S_{[n]\setminus\Delta}$ and $S_{\Theta}\times S_{[n]\setminus\Theta}$ containing elements with cycle-type $(n-3,2,1)$ and the fixed point $x$ are in $Max(S_n,H)$, such that $\Delta\cap\Psi=\{x,y\}$ and $\Theta\cap\Psi=\{x,z\}$. Hence, $x$ is the fixed point of $(S_{\Delta}\times S_{[n]\setminus\Delta})\cap(S_{\Theta}\times S_{[n]\setminus\Theta})$ and thus $x$ is the fixed point of $H$. This leads to our observation. The proof of this lemma is complete.   \qed

Using Lemma\ \ref{pan2-1} and Lemma\ \ref{pan2-2}, we obtain the following theorem immediately.

\begin{thm}\label{pan2-3}\normalfont
All nontrivial subgroups of $S_n$ are FF-subgroups.
\end{thm}

We are now turning to study the alternating groups. The proof of alternating group is similar to that given earlier for symmetric group and so we will not go into details here.

\begin{thm}\label{pan2-4}\normalfont
All nontrivial subgroups of $A_n$ are FF-subgroups.
\end{thm}
\demo Assume that the nontrivial subgroup $H$ of $A_n$ is not a FF-subgroup. Case 1: $n$ is an odd number. Proceeding as in the proof of Lemma\ \ref{pan2-2}, we see all elements with cycle-type $(n-3,2,1)$ are contained in the maximal subgroups $A_{n-1}$ or $A_3\times A_{n-3}$ or $A_2\times A_{n-2}$. Using the same argument as in the proof of Lemma\ \ref{pan2-2} to show that $H$ is the identity. Case 2: $n$ is an even number. According to the proof of Lemma\ \ref{pan2-1}, it follows that all elements with cycle-type $(n-2,2)$ are contained in the maximal subgroups $A_2\times A_{n-2}$. Similarly, we can obtain that $H$ is the identity. The proof of this theorem is complete.  \qed

\section {Projective Groups}

Recall that the projective group $PSL(d,q)$ has a faithful 2-transitive action of degree $\frac{q^d-1}{q-1}$ on the set of 1-dimensional subspaces of $F_q$. However, the maximal subgroups of $PSL(d,q)$ are rather complicated except $d=2$, see the references \cite{A,KL}. So we mainly investigate $PSL(2,q)$ in this section. See the reference \cite{D}, the maximal subgroups of $PSL(2,q)$ are as follows.

\begin{lem}\label{pan3-1}\normalfont(\cite{D}~)
(A). Let $q=2^f\geq4$. Then the maximal subgroup of $PSL(2,q)$ is one of the following:\\
1) $Z^f_2\rtimes Z_{q-1}$, that is, the stabilizer subgroup of the projective line $PG(1,q)$;\\
2) $D_{2(q-1)}$;\\
3) $D_{2(q+1)}$;\\
4) $PGL(2,q_0)$, $q=q^r_0$ and $q_0\neq2$ where $r$ is a prime number.\\
(B). Let $q=p^f\geq5$, where $p$ is an odd prime number. Then the maximal subgroup of $PSL(2,q)$ is one of the following:\\
1) $Z^f_p\rtimes Z_{(q-1)/2}$, that is, the stabilizer subgroup of the projective line $PG(1,q)$;\\
2) $D_{q-1}$ with $q\geq13$;\\
3) $D_{q+1}$, $q\neq7,9$;\\
4) $PGL(2,q_0)$, $q=q^2_0$;\\
5) $PSL(2,q_0)$, $q=q^r_0$ where $r$ is an odd prime number;\\
6) $A_5$, $q=p$, $q\equiv\pm1(mod~10)$ or $q=p^2$, $p\equiv\pm3(mod~10)$;\\
7) $A_4$, $q=p\equiv\pm3(mod~8)$ and $q\equiv\pm1(mod~10)$;\\
8) $S_4$, $q=p\equiv\pm1(mod~8)$.
\end{lem}

Now we start to study the FF-subgroups of $PSL(2,q)$ and further give the following lemmas.
\begin{lem}\label{pan3-2}\normalfont
Let $q=2^f\geq4$. Then all nontrivial subgroups of $PSL(2,q)$ are FF-subgroups.
\end{lem}
\demo We claim that all elements of order $q+1$ are contained in the maximal subgroups $D_{2(q+1)}$. It follows from Lemma\ \ref{pan3-1} (A) that there are four classes of maximal subgroups in $PSL(2,q)$, and then we verify our claim one by one.

Obviously, there exist some elements of order $q+1$ in the maximal subgroup $D_{2(q+1)}$. Consider other classes of maximal subgroups. Since $q=2^f\geq4$, it follows that $q+1$ and $q-1$ are coprime, and thus there is no element of order $q+1$ in the maximal subgroups $Z^f_2\rtimes Z_{q-1}$ and $D_{2(q-1)}$. It is well-known that $|PGL(2,q_0)|=\frac{(q^2_0-1)(q^2_0-q_0)}{q_0-1}=q_0(q^2_0-1)$. However, $q=q^r_0$ and $q_0\neq2$ indicate $q+1\nmid |PGL(2,q_0)|$, and so there is no element of order $q+1$ in $PGL(2,q_0)$. So far, we have obtained our claim.

Suppose the nontrivial subgroup $H$ of $PSL(2,q)$ is not a FF-subgroup. Thus $\Delta_H(PSL(2,q))=PSL(2,q)$. Then by our claim and \cite[Theorem~8.3,~8.5]{H} we see that all maximal subgroups of the type $D_{2(q+1)}$ are in $Max(PSL(2,q),H)$ and further $H$ is the identity, a contradiction. This leads to the lemma. \qed

\begin{lem}\label{pan3-3}\normalfont
Let $q=p^f\geq5$, where $p$ is an odd prime number. Then all nontrivial subgroups of $PSL(2,q)$ are FF-subgroups.
\end{lem}
\demo Consider the case that $q<13$. In this case, we have $q=5,7,9,11$. According to \cite[Chaper~2,~Theorem~6.14]{H} and Theorem\ \ref{pan2-4}, it follows that $PSL(2,9)\cong A_6$ and $PSL(2,5)\cong A_5$, and further the lemma holds for $q=9$ and $q=5$. If $q=p=7$, then the maximal subgroup is 1) or 8) of Lemma\ \ref{pan3-1} (B), and further all elements of order $4$ are contained in the maximal subgroup 8); If $q=p=11$, then the maximal subgroup is 1) or 3) or 6) or 7) of Lemma\ \ref{pan3-1} (B), and further all elements of order $6$ are contained in the maximal subgroup 3). Then by using the similar argument as in proof of Lemma\ \ref{pan3-2}, we see the lemma is also true for $q=p=7$ and $q=p=11$.

In the case of $q\geq13$, we claim that all elements of order $\frac{q+1}{2}$ are contained in the maximal subgroup $D_{(q+1)}$. One easily checks that there is no element of order $\frac{q+1}{2}$ in the maximal subgroup 1), 2), 6), 7), 8) of Lemma\ \ref{pan3-1} (B). Moreover, it is well-known that $|PGL(2,q_0)|=\frac{(q^2_0-1)(q^2_0-q_0)}{q_0-1}=q_0(q^2_0-1)$ and $|PSL(2,q_0)|=\frac{|PGL(2,q_0)|}{2}$, and so $\frac{q+1}{2}\nmid |PGL(2,q_0)|$ and $\frac{q+1}{2}\nmid |PSL(2,q_0)|$. Hence, our claim holds. An argument similar to the one used in the proof of Lemma \ \ref{pan3-2} shows this lemma is true for $q\geq13$. The proof of this lemma is complete.    \qed

Applying the Lemma \ \ref{pan3-2} and Lemma \ \ref{pan3-3}, we obtain the following theorem.

\begin{thm}\label{pan3-4}\normalfont
All nontrivial subgroups of $PSL(2,q)$ are FF-subgroups.
\end{thm}

\section {Final remarks}

Note that for the $2$-generator group, the probability of its nontrivial subgroups are FF-subgroups is very high. On the other hand, the Frattini subgroup of finite simple group is trivial. So we propose the following conjecture and open problem to finish this paper.

\begin{Con}\label{pan4-1}\normalfont
All nontrivial subgroups of a finite simple group are FF-subgroups.
\end{Con}

\begin{Que}\label{pan4-2}\normalfont
Is there a $2$-generator group has a nontrivial subgroup is neither the subgroup of its Frattini subgroup nor a FF-subgroup.
\end{Que}

\end{document}